\documentclass{article}
\usepackage{amssymb}
\usepackage{amsmath}
\usepackage{color}
\newtheorem{theorem}{Theorem}

\newtheorem{proposition}[theorem]{Proposition}
\newtheorem{remark}[theorem]{Remark}

\newenvironment{proof}[1][Proof]{\noindent\textbf{#1.} }{\ \rule{0.5em}{0.5em}}
\begin{document}

\title{A Note on Carlier's Inequality}
\author{\textbf{Regina S. Burachik} \\
\textit{Mathematics, UniSA STEM}\\
\textit{University of South Australia}\\
\textit{Mawson Lakes, SA 5095, Australia}\\
\textit{Regina.Burachik@unisa.edu.au} \and \textbf{J. E. Mart\'{\i}nez-Legaz}
\\
\textit{Department of Economics and Economic History}\\
\textit{Autonomous University of Barcelona}\\
\textit{Bellaterra, Spain}\\
\textit{juanenrique.martinez.legaz@uab.cat}}
\date{}
\maketitle

\begin{abstract}
Recently, Carlier established in \cite{Carlier} a quantitave version of
Fitzpatrick inequality in a Hilbert space. We extend Carlier's result to the
framework of reflexive Banach spaces. In the Hilbert space setting, we
obtain an improved version of the strong Fitzpatrick inequality due to
Voisei and Z\u{a}linescu.
\end{abstract}

\section{Introduction}

Let $H$ be a Hilbert space with duality pairing denoted by $\langle \cdot
,\cdot \rangle $ and norm denoted by $\Vert \cdot \Vert $. For a maximally
monotone map $T:H\rightrightarrows H$, Carlier proved in \cite[eq 2.3]%
{Carlier} the following inequality for every $\lambda>0$: 
\begin{equation}
\mathcal{F}_{T}(x,v)-\left\langle x,v\right\rangle \geq \dfrac{1}{\lambda }%
\Vert x-(I+\lambda T)^{-1}(x+\lambda v)\Vert ^{2},  \label{eq1}
\end{equation}%
where $\mathcal{F}_{T}:H^{2}\rightarrow \mathbb{R\cup }\left\{ +\infty
\right\} $ is the \emph{Fitzpatrick function} \cite{F88} associated with $T,$
defined by%
\begin{equation*}
\mathcal{F}_{T}(x,v):=\langle x,v\rangle -\inf_{\left( z,w\right) \in
G\left( T\right) }\langle x-z,v-w\rangle ;
\end{equation*}%
here we use $G\left( T\right) $ to denote the graph of $T,$ given by%
\begin{equation*}
G\left( T\right) :=\left\{ \left( z,w\right) \in H^{2}:w\in T\left( z\right)
\right\} .
\end{equation*}

Inequality \eqref{eq1} constitutes a strengthening of \emph{Fitzpatrick's
inequality}%
\begin{equation*}
\mathcal{F}_{T}(x,v)-\left\langle x,v\right\rangle \geq 0.
\end{equation*}%
The aim of this paper is twofold:\ to extend \emph{Carlier's inequality}
beyond the Hilbert space setting and, using the same technique, to improve,
in the framework of Hilbert spaces, the strong Fitzpatrick inequality due to
Voisei and Z\u{a}linescu \cite[Theorem 2.6]{VZ09}.

Recall that, for a general reflexive Banach space $\left( X,\left\Vert \cdot
\right\Vert \right) $ with dual $\left( X^{\ast },\left\Vert \cdot
\right\Vert _{\ast }\right) $ and duality product $\langle \cdot ,\cdot
\rangle :X\times X^{\ast }\rightarrow \mathbb{R},$ the \emph{normalized
duality mapping} is defined as the point to set operator $%
J_{X}:X\rightrightarrows X^{\ast }$ given by 
\begin{equation*}
J_{X}(x):=\{v\in X^{\ast }:\left\langle x,v\right\rangle =\Vert x\Vert \Vert
v\Vert _{\ast },\,\Vert x\Vert =\Vert v\Vert _{\ast }\}.
\end{equation*}%
This normalized duality mapping plays the role of the identity map. Indeed, $%
J_{X}$ is the subdifferential operator of the function $\frac{1}{2}\Vert
\cdot \Vert ^{2}$. We recall that the \emph{subdifferential} operator $%
\partial f$ of a function $f:X\rightarrow \mathbb{R}\mathbb{\cup }\left\{
+\infty \right\} $ is defined by%
\begin{equation*}
\partial f\left( x\right) :=\left\{ v\in X^{\ast }:f\left( y\right) \geq
f\left( x\right) +\left\langle y-x,v\right\rangle \text{ }\forall y\in
X\right\} .
\end{equation*}%
It is well known that $\partial f$ is maximally monotone and has the
function $h:X\times X^{\ast }\rightarrow \mathbb{R\cup }\left\{ +\infty
\right\} $ defined by $h\left( x,v\right) :=f\left( x\right) +f^{\ast }\left(
v\right) $ as a \emph{convex representation}, that is, it is convex and norm$%
\times $weak$^{\ast }$ lower semicontinuous and satisfies the following two
conditions:

\begin{enumerate}
\item $h(x,v)\geq \left\langle x,v\right\rangle ,$\quad $\forall (x,v)\in
H^{2},$

\item $h(x,v)=\left\langle x,v\right\rangle \Leftrightarrow (x,v)\in G\left(
T\right) .$
\end{enumerate}

In the definition of $h$ above, $f^{\ast }:X^{\ast }\rightarrow \mathbb{%
R\cup }\left\{ +\infty \right\} $ denotes the \emph{Fenchel conjugate}
function of $f,$ given by%
\begin{equation*}
f^{\ast }\left( v\right) :=\sup_{x\in X}\left\{ \left\langle
x,v\right\rangle -f\left( x\right) \right\} ,
\end{equation*}

It can be proved directly from the definitions that the Fitzpatrick function $\mathcal{F}_{T}$ is finite everywhere when $T:=J_X$, i.e., when $T$
is the normalized duality mapping. In a Hilbert space $H$,
since the dual space $H^{\ast }$ is canonically identified with $H,$ the
normalized duality mapping is precisely the identity map. In such a setting,
the normalized duality mapping is strongly monotone. Indeed, recall that a
maximally monotone operator $B:X\rightrightarrows X^{\ast }$ is \emph{%
strongly monotone} with constant $c>0$ if, for every $v\in Bx$, $u\in By$,
we have 
\begin{equation}
\left\langle x-y,v-u\right\rangle \geq c\Vert x-y\Vert ^{2}.  \label{str mon}
\end{equation}%
It is straightforward to check that the identity map of a Hilbert space is
strongly monotone with constant $c=1$.

It turns out that these two properties, namely the finiteness of the
Fitzpatrick function and the strong monotonicity of a maximally monotone
operator, are what we need in order to achieve an extension of Carlier's
inequality to a reflexive Banach space. Namely, if a maximally monotone
operator $B$ satisfies these two properties, then we can use $B$ in place of
the identity map to obtain a more general version of (\ref{eq1}). In doing
so, we are inspired by the analysis in \cite{M08}. In the case of a $2$%
-uniformly convex Banach space, one can take $B:=J_{X},$ the normalized
duality mapping.

We will also show that a strengthened version of the strong Fitzpatrick
inequality mentioned above can be easily obtained, in the Hilbert space
framework, using the same technique as in the proof of Carlier's inequality.

The present note is organized as follows. In the next section, we extend
inequality (\ref{eq1}) to the context of reflexive Banach spaces. In Section %
\ref{2-unif}, we prove that the Fenchel subdifferential operator of a
strongly convex function on a $2$-uniformly convex Banach space is strongly
monotone and use this result to express, for maximally monotone operators
defined on such spaces, our Carlier type inequality in terms of the
normalized duality mapping. Finally, in the last section we obtain an
improved version of the strong Fitzpatrick inequality.

\section{The case of a reflexive Banach space\label{Sect refl}}

In this section we assume that $\left( X,\left\Vert \cdot \right\Vert
\right) $ is a reflexive real Banach space. Let $T,B:X\rightrightarrows
X^{\ast }$ be maximally monotone, and assume that $B$ has a finite-valued
Fitzpatrick function and is strongly monotone with constant $c.$ By \cite[%
Corollary 2.7]{M08}, for $\lambda >0$ the operator $\left( B+\lambda
T\right) ^{-1}$ has full domain;\ moreover, it is single valued because $B$
is strictly monotone. Let $\left( x,v\right) \in X\times X^{\ast }.$ Given $%
x^{\prime }\in Bx$, we define $w_{\lambda }:=\left( B+\lambda T\right)
^{-1}\left( x^{\prime }+\lambda v\right) .$ We then have $x^{\prime
}+\lambda v\in Bw_{\lambda }+\lambda Tw_{\lambda },$ that is, $x^{\prime
}+\lambda v\in w_{\lambda }^{\prime }+\lambda Tw_{\lambda }$ for some $%
w_{\lambda }^{\prime }\in Bw_{\lambda }.$ It follows that $\frac{x^{\prime
}-w_{\lambda }^{\prime }}{\lambda }+v\in Tw_{\lambda },$ which means that $%
\left( w_{\lambda },\frac{x^{\prime }-w_{\lambda }^{\prime }}{\lambda }%
+v\right) \in G\left( T\right) .$ Therefore, we have 
\begin{eqnarray}
F_{T}\left( x,v\right) -\left\langle x,v\right\rangle  &=&\sup_{\left(
y,w\right) \in G\left( T\right) }\left\langle x-y,w-v\right\rangle 
\label{ineq} \\
&\geq &\left\langle x-w_{\lambda },\frac{x^{\prime }-w_{\lambda }^{\prime }}{%
\lambda }+v-v\right\rangle   \notag \\
&=&\frac{1}{\lambda }\left\langle x-w_{\lambda },x^{\prime }-w_{\lambda
}^{\prime }\right\rangle\\
& \geq& \frac{c}{\lambda }\left\Vert x-w_{\lambda
}\right\Vert ^{2}  \notag
=\frac{c}{\lambda }\left\Vert x-\left( B+\lambda T\right) ^{-1}\left(
x^{\prime }+\lambda v\right) \right\Vert ^{2},  \notag
\end{eqnarray}%
where the second inequality follows from (\ref{str mon}). By taking the supremum
over $x^{\prime }\in Bx,$ we conclude that%
\begin{equation*}
F_{T}\left( x,v\right) -\left\langle x,v\right\rangle \geq \frac{c}{\lambda }%
\sup_{x^{\prime }\in Bx}\left\Vert x-\left( B+\lambda T\right) ^{-1}\left(
x^{\prime }+\lambda v\right) \right\Vert ^{2}.
\end{equation*}%

We have thus proved the following theorem.

\begin{theorem}
\label{GCI}Let $T,B:X\rightrightarrows X^{\ast }$ be maximally monotone.
Assume that $B$ is strongly monotone with constant $c$ and has a
finite-valued Fitzpatrick function $F_{T}$. Then, for every $\left(
x,v\right) \in X\times X^{\ast }$ and $\lambda >0,$ one has%
\begin{equation*}
F_{T}\left( x,v\right) -\left\langle x,v\right\rangle \geq \frac{c}{\lambda }%
\sup_{x^{\prime }\in Bx}\left\Vert x-\left( B+\lambda T\right) ^{-1}\left(
x^{\prime }+\lambda v\right) \right\Vert ^{2}.
\end{equation*}
\end{theorem}

Clearly, Carlier's inequality immediately follows from Theorem \ref{GCI},
since, when $X$ is a Hilbert space, using the canonical identification of $%
X^{\ast }$ with $X$ one can take $B$ to be the identity mapping $I$, which
is strongly monotone with constant $c=1$. Indeed, by taking $B:=I$ (the single valued identity map) we readily obtain \eqref{eq1}. Recall also that the identity map's Fitzpatrick function is given by%
\begin{equation*}
F_{I}\left( x,v\right) =\frac{1}{4}\left\Vert x+v\right\Vert ^{2}.
\end{equation*}

\section{The case of a $2$-uniformly convex Banach space\label{2-unif}}

We recall that a Banach space $\left( X,\left\Vert \cdot \right\Vert \right) 
$ is said to be $2$\emph{-uniformly convex} with constant $\mu >0$ if the
following implication holds true:%
\begin{equation*}
x,y\in X,\text{ }\left\Vert x\right\Vert =1=\left\Vert y\right\Vert
\Rightarrow \left\Vert \frac{x+y}{2}\right\Vert \leq 1-\mu \left\Vert
x-y\right\Vert ^{2}.
\end{equation*}%
This is equivalent to saying that the \emph{modulus of convexity }$\delta
_{X}:\left[ 0,2\right] \rightarrow \mathbb{R}$ of $X$, defined by%
\begin{equation*}
\delta _{X}\left( \epsilon \right) :=\inf \left\{ 1-\left\Vert \frac{x+y}{2}%
\right\Vert :x,y\in X,\text{ }\left\Vert x\right\Vert =1=\left\Vert
y\right\Vert ,\text{ }\left\Vert x-y\right\Vert \geq \epsilon \right\} ,
\end{equation*}%
satisfies the inequality%
\begin{equation}
\delta _{X}\left( \epsilon \right) \geq \mu \epsilon ^{2}  \label{2uc}
\end{equation}%
for every $\epsilon \in \left[ 0,2\right] .$

It is well known and easy to see, using the parallelogram law, that every
Hilbert space is $2$-uniformly convex with constant $\frac{1}{8}.$ Other
important examples of $2$-uniformly convex spaces are the $L^{p}$ spaces for 
$p\in \left( 1,2\right] .$ Indeed, for such spaces (\ref{2uc}) holds with $%
\mu =\frac{p-1}{8};$ this can be easily proved using \cite[Proposition 3]%
{BCL94}.

To consider the special case of Theorem \ref{GCI} in which $B$ is a
subdifferential operator, we need the following preliminary result on
strongly convex functions. Recall that $f:X\rightarrow \mathbb{R}$ is \emph{%
strongly convex} with constant $m>0$ if $f-\frac{m}{2}\left\Vert \cdot
\right\Vert ^{2}$ is convex. 

The proof of the following result is inspired by that of \cite[Proposition
2.11]{PR86};\ in particular, we will use the fact, proved there, that, if $%
\left( X,\left\Vert \cdot \right\Vert \right) $ is $2$\emph{-uniformly convex%
} with constant $\mu ,$ then, for every $x,y\in X,$ one has%
\begin{equation*}
\left\Vert x\right\Vert ^{2}+\left\Vert y\right\Vert ^{2}-\frac{1}{2}%
\left\Vert x+y\right\Vert ^{2}\geq \frac{\mu }{2}\left\Vert x-y\right\Vert
^{2},
\end{equation*}%
an inequality which follows easily from \cite[Lemma 1.e.10]{LT73}.

\begin{proposition}
\label{str conv}Let $\left( X,\left\Vert \cdot \right\Vert \right) $ be $2$%
-uniformly convex with constant $\mu >0$ and $f:X\rightarrow \mathbb{R}$ be
strongly convex with constant $m>0$. 
Then 
$\partial f$ is strongly monotone with constant $\frac{m\mu }{2}$.
\end{proposition}

\begin{proof}
Let $x,y\in X,$ $v\in \partial f\left( x\right) $ and $w\in \partial f\left(
y\right) .$ From the subgradient inequalities $f\left( \frac{x+y}{2}\right)
\geq f\left( x\right) +\frac{1}{2}\left\langle y-x,v\right\rangle $ and $%
f\left( \frac{x+y}{2}\right) \geq f\left( y\right) +\frac{1}{2}\left\langle
x-y,w\right\rangle ,$ using the strong convexity of $f$ we obtain%
\begin{eqnarray*}
\left\langle x-y,v-w\right\rangle &\geq &4\left( \frac{f\left( x\right)
+f\left( y\right) }{2}-f\left( \frac{x+y}{2}\right) \right) \\
&\geq &m\left( \left\Vert x\right\Vert ^{2}+\left\Vert y\right\Vert ^{2}-%
\frac{1}{2}\left\Vert x+y\right\Vert ^{2}\right) \geq \frac{m\mu }{2}%
\left\Vert x-y\right\Vert ^{2}.
\end{eqnarray*}
\end{proof}

\bigskip

The normalized duality mapping $J_{X}$ is the subdifferential operator of $%
\frac{1}{2}\left\Vert \cdot \right\Vert ^{2},$ which is a strongly convex
function with constant $1$. Moreover, $F_{J_{X}}$ is finite-valued, since,
for every $\left( x,v\right) \in X\times X^{\ast },$ one has $%
F_{J_{X}}\left( x,v\right) \leq \frac{1}{2}\left( \left\Vert x\right\Vert
^{2}+\left\Vert v\right\Vert _{\ast }^{2}\right) ,$ as the right hand side
of this inequality defines a convex representation of $J_{X}$ and $F_{J_{X}}$
is the smallest of such representations. On the other hand, the Milman--Pettis
theorem states that every uniformly convex Banach space is reflexive.
Therefore, Theorem \ref{GCI} and Proposition \ref{str conv} yield the
following result.

\begin{theorem}
Let $\left( X,\left\Vert \cdot \right\Vert \right) $ be $2$-uniformly convex
with constant $\mu >0$ and $T:X\rightrightarrows X^{\ast }$ be maximally
monotone. Then, for every $\left( x,v\right) \in X\times X^{\ast },$ one has%
\begin{equation*}
F_{T}\left( x,v\right) -\left\langle x,v\right\rangle \geq \frac{\mu }{%
2\lambda }\sup_{x^{\prime }\in J_{X}(x)}\left\Vert x-\left( J_{X}+\lambda
T\right) ^{-1}\left( x^{\prime }+\lambda v\right) \right\Vert ^{2}.
\end{equation*}
\end{theorem}

\section{An improved version of the strong Fitzpatrick inequality in Hilbert
spaces\label{Fitz}}

Another strengthening of Fitzpatrick inequality, called the strong
Fitzpatrick inequality, was obtained in 2009 by Voisei and Z\u{a}linescu 
\cite[Theorem 2.6]{VZ09} (see also \cite[Theorem 9.7.2]{BV10}) for maximally
monotone operators of type (NI) (see, e.g., \cite{S98} for the definition of
this class of operators). In the case of reflexive spaces, every maximally
monotone operator is of type (NI)\ \cite[p. 149]{S98}; for such spaces, a
more general version of the strong Fitzpatrick inequality was given in \cite[%
Theorem 5]{Elias-Martinez Legaz-2017}. Unlike Carlier's inequality, the
strong Fitzpatrick inequality is expressed in terms of the distance to the
graph of the given operator. In this section, using the same argument as in
the proof of Theorem \ref{GCI} we show how to derive a stronger version of
the strong Fitzpatrick inequality for maximally monotone operators defined
on a Hilbert space.

In what follows, we assume that $H$ is a Hilbert space with duality pairing
denoted by $\langle \cdot ,\cdot \rangle $ and norm denoted by $\Vert \cdot
\Vert $. We use (\ref{ineq}) with $B$ equal to the identity mapping, which
is obviously strongly monotone with constant $c=1$. Consider $w_{\lambda}$ as in the beginning of Section \ref{Sect refl}. From the inequality $%
F_{T}\left( x,v\right) -\left\langle x,v\right\rangle \geq \frac{1}{\lambda }%
\left\Vert x-w_{\lambda }\right\Vert ^{2}$, one obtains%
\begin{eqnarray*}
F_{T}\left( x,v\right) -\left\langle x,v\right\rangle &\geq &\frac{1}{%
\lambda }\left( \frac{\lambda ^{2}}{1+\lambda ^{2}}\left\Vert x-w_{\lambda
}\right\Vert ^{2}+\frac{1}{1+\lambda ^{2}}\left\Vert x-w_{\lambda
}\right\Vert ^{2}\right) \\
&=&\frac{\lambda }{1+\lambda ^{2}}\left( \left\Vert x-w_{\lambda
}\right\Vert ^{2}+\left\Vert \frac{x-w_{\lambda }}{\lambda }\right\Vert
^{2}\right) \\
&=&\frac{\lambda }{1+\lambda ^{2}}\left( \left\Vert x-w_{\lambda
}\right\Vert ^{2}+\left\Vert \frac{x-w_{\lambda }}{\lambda }+v-v\right\Vert
^{2}\right) \\
&\geq &\frac{\lambda }{1+\lambda ^{2}}\inf_{\left( w,z\right) \in G\left(
T\right) }\left\{ \left\Vert x-w\right\Vert ^{2}+\left\Vert v-z\right\Vert
^{2}\right\},
\end{eqnarray*}
where in the last inequality we used again (as in Section \ref{Sect refl}), the fact that $\dfrac{x-w_{\lambda}}{\lambda}+v\in Tw_{\lambda}$. Setting $\lambda :={1},$ the following result is proved.

\begin{theorem}
Let $H$ be a Hilbert space. If $T:H\rightrightarrows H$ is maximally
monotone, then, for every $\left( x,v\right) \in H\times H,$ one has%
\begin{equation}
F_{T}\left( x,v\right) -\left\langle x,v\right\rangle \geq \frac{1}{2}%
\inf_{\left( w,z\right) \in G\left( T\right) }\left\{ \left\Vert
x-w\right\Vert ^{2}+\left\Vert v-z\right\Vert ^{2}\right\} .  \label{SFI}
\end{equation}
\end{theorem}

\begin{remark}
Inequality (\ref{SFI}) is stronger than the strong Fitzpatrick inequality
obtained by Voisei and Z\u{a}linescu \cite[Theorem 4]{VZ09} (see \cite[%
Theorem 9.7.2]{BV10}), as in the latter the coefficient $\frac{1}{4}$
appears in place of our $\frac{1}{2}.$
\end{remark}

\textbf{Acknowledgment.} Juan Enrique Mart\'{\i}nez-Legaz has been partially
supported by Grant PID2022-136399NB-C22 from MICINN, Spain, and ERDF,
\textquotedblright A way to make Europe", European Union.

\end{document}